\documentclass[11pt]{amsart}

\usepackage{texdraw,multicol,rotating}
\usepackage{amsmath}
\usepackage{amssymb}
\usepackage{epsfig}
\usepackage{amssymb,amsthm,amsfonts,amscd}
\input{txdtools.tex}
\theoremstyle{plain}
\newtheorem{thm}{Theorem}[section]
\newtheorem{prop}[thm]{Proposition}

\newtheorem{cor}[thm]{Corollary}
\theoremstyle{definition}
\newtheorem{defi}[thm]{Definition}

\newtheorem{example}[thm]{Example}
\theoremstyle{remark}
\newtheorem{remark}[thm]{Remark}
\numberwithin{equation}{section}
\newbox{\tmpa}
\newbox{\tmpb}

\newcommand{\nc}{\newcommand}
\nc{\Uq}{U_q} \nc{\Z}{\mathbf{Z}} \nc{\C}{\mathbf{C}}
\nc{\Q}{\mathbf{Q}}
\nc{\op}{\oplus} \nc{\ot}{\otimes} \nc{\pv}{P^{\vee}}
\nc{\ali}{\alpha_i} \nc{\B}{\mathbf{B}} \nc{\F}{\mathbf{F}}
\nc{\bP}{\mathbf{P}} \nc{\V}{\mathbf{V}} \nc{\La}{\Lambda}
\nc{\la}{\lambda}
\nc{\nbinom}[2]{\genfrac{}{}{0pt}{1}{{#1}}{{#2}}}
\nc{\qbinom}[2]{\left[\genfrac{}{}{0pt}{1}{{#1}}{{#2}}\right]}
\nc{\path}{\mathcal{P}} \nc{\fit}{\tilde{f}_i}
\nc{\eit}{\tilde{e}_i} \nc{\Y}{\mathbf{Y}} \nc{\A}{\mathbf{A}}
\nc{\ra}{\rightarrow} \nc{\vep}{\varepsilon} \nc{\vphi}{\varphi}
\nc{\g}{\mathfrak{g}} \nc{\h}{\mathfrak{h}} \nc{\oP}{\overline{P}}
\nc{\pathp}{\mathbf{p}}
\nc{\tris}{ \bsegment \move(0 0)\lvec(10 0)\lvec(10 10)\lvec(0
0)\ifill f:0.7 \esegment } \nc{\recs}{ \bsegment \move(0
0)\lvec(10 0)\lvec(10 5)\lvec(0 5)\lvec(0 0)\ifill f:0.7 \esegment
}
\nc{\hcvec}[5]{%
\getpos(#1 #3)\spx\spy \getpos(#2 #3)\epx\epy \getpos(#4
#5)\xoff\yoff \realadd \spx \xoff \twox \realadd \epx {-\xoff}
\thrx \realadd \spy \yoff \posy \move({\spx} {\spy}) \clvec
({\twox} {\posy})({\thrx} {\posy})({\epx} {\epy}) \rmove(0 0) }
\nc{\ahead}[2]{%
\cossin (0 0)({#1} {#2})\cosa\sina \bsegment
  \drawdim in \setunitscale 0.065
  \realmult {-0.5} \cosa \hcosa
  \realmult {-0.5} \sina \hsina
  \move({\hcosa} {\hsina}) \ravec({\cosa} {\sina})
\esegment }
\nc{\boxi}{%
{%
\savebox{\tmppic}{\begin{texdraw} \small \drawdim em \textref h:C
v:C \setunitscale 0.55 \htext(0 0){$i$} \move(-1 -1)\lvec(-1
1)\lvec(1 1)\lvec(1 -1)\lvec(-1 -1)
\end{texdraw}}%
\raisebox{-0.19\height}{\usebox{\tmppic}}%
}%
}
\nc{\boxj}{%
{%
\savebox{\tmppic}{\begin{texdraw} \small \drawdim em \textref h:C
v:C \setunitscale 0.55 \htext(0 0.1){$j$} \move(-1 -1)\lvec(-1
1)\lvec(1 1)\lvec(1 -1)\lvec(-1 -1)
\end{texdraw}}%
\raisebox{-0.19\height}{\usebox{\tmppic}}%
}%
}
\nc{\boxipo}{%
{%
\savebox{\tmppic}{\begin{texdraw} \small \drawdim em \textref h:C
v:C \setunitscale 0.55 \htext(0.15 0){$i\!\!+\!\!1$} \move(-1.4
-1)\lvec(-1.4 1)\lvec(1.4 1)\lvec(1.4 -1)\lvec(-1.4 -1)
\end{texdraw}}%
\raisebox{-0.19\height}{\usebox{\tmppic}}%
}%
} \everytexdraw{ \drawdim in \arrowheadsize l:0.065 w:0.03
\arrowheadtype t:F \textref h:C v:C }
\newsavebox{\tmppic}
\newsavebox{\tmpfig}
\newsavebox{\tmpdraw}
\newsavebox{\tmpfiga}
\newsavebox{\tmpfigb}
\newsavebox{\tmpfigc}
\newsavebox{\tmpfigd}
\newsavebox{\tmpfige}
\newsavebox{\tmpfigf}
\newsavebox{\tmpfigg}
\newsavebox{\tmpfigh}
\newsavebox{\tmpfigi}
\newsavebox{\tmpfigj}
\newsavebox{\tmpfigk}
\newsavebox{\tmpfigl}
\newsavebox{\tmpfigm}
\newsavebox{\tmpfign}
\newsavebox{\tmpfigo}
\newsavebox{\tmpfigp}
\newsavebox{\tmpfigq}
\newsavebox{\tmpfigr}
\newsavebox{\tmpfigs}
\newsavebox{\tmpfigt}
\newsavebox{\tmpfigu}
\newsavebox{\tmpfigv}
\newsavebox{\tmpfigw}
\newsavebox{\tmpfigx}
\newsavebox{\tmpfigy}
\newsavebox{\tmpfigz}
\newsavebox{\tmpfigaa}
\newsavebox{\tmpfigab}
\newsavebox{\tmpfigac}
\newsavebox{\tmpfigad}
\newsavebox{\tmpfigae}
\newsavebox{\tmpfigaf}
\newsavebox{\tmpfigag}
\newsavebox{\tmpfigah}
\newsavebox{\tmpfigai}
\newsavebox{\tmpfigaj}
\newsavebox{\tmpfigak}
\newsavebox{\tmpfigal}
\newsavebox{\tmpfigam}
\newsavebox{\tmpfigan}
\newsavebox{\tmpfigao}
\newsavebox{\tmpfigap}
\newsavebox{\tmpfigaq}
\newsavebox{\tmpfigar}
\newsavebox{\tmpfigas}
\newsavebox{\tmpfigat}
\newsavebox{\tmpfigau}
\newsavebox{\tmpfigav}
\newsavebox{\tmpfigaw}
\newsavebox{\tmpfigax}
\newsavebox{\tmpfigay}
\newsavebox{\tmpfigaz}
\newsavebox{\tmpfigba}
\newsavebox{\tmpfigbb}
\newsavebox{\tmpfigbc}
\newsavebox{\tmpfigbd}
\newsavebox{\tmpfigbe}
\newsavebox{\tmpfigbf}
\newsavebox{\tmpfigbg}
\newsavebox{\tmpfigbh}

\newsavebox{\spinaa}
\newsavebox{\spinab}
\newsavebox{\spinac}
\newsavebox{\spinad}
\newsavebox{\spinae}
\newsavebox{\spinaf}
\newsavebox{\spinag}
\newsavebox{\spinah}
\newsavebox{\spinai}
\newsavebox{\spinaj}
\newsavebox{\spinak}
\newsavebox{\spinal}
\newsavebox{\spinam}
\newsavebox{\spinba}
\newsavebox{\spinbb}
\newsavebox{\spinbc}
\newsavebox{\spinbd}
\newsavebox{\spinbe}
\newsavebox{\spinbf}
\newsavebox{\spinbg}
\newsavebox{\spinbh}
\nc{\node}{\lcir r:1 }
\nc{\sline}{\bsegment\savepos(10 0)(*ex *ey)
            \move(1 0)\rlvec(8 0)
            \esegment\move(*ex *ey)}
\nc{\dline}{\bsegment\savepos(10 0)(*ex *ey)
            \move(0.93 0.4)\rlvec(8.14 0)\rmove(0 -0.8)\rlvec(-8.14 0)
            \esegment\move(*ex *ey)}
\nc{\uline}{\bsegment\savepos(0 10)(*ex *ey)
            \move(0 1)\rlvec(0 8)
            \esegment\move(*ex *ey)}
\nc{\lpoint}{\savecurrpos(*ex *ey)
             \rmove(2.5 1.5)\rlvec(-1.5 -1.5)\rlvec(1.5 -1.5)
             \move(*ex *ey)}
\nc{\rpoint}{\savecurrpos(*ex *ey)
             \rmove(-2.5 -1.5)\rlvec(1.5 1.5)\rlvec(-1.5 1.5)
             \move(*ex *ey)}
\nc{\bline}{\bsegment\savepos(10 0)(*ex *ey)
            \linewd 0.6 \move(1.1 0)\rlvec(7.8 0)
            \esegment\move(*ex *ey)}
\nc{\araise}[1]{\raisebox{4.5pt}{#1}}
\nc{\braise}[1]{\raisebox{12.1pt}{#1}}
\nc{\craise}[1]{\raisebox{8pt}{#1}}
\nc{\draise}[1]{\raisebox{12pt}{#1}}
\nc{\eraise}[1]{\raisebox{14.5pt}{#1}}

\begin{document}

\title[Monomial Relization of Crystal Bases for Special Linear Lie
Algebras]
{Monomial Relization of Crystal Bases for Special Linear
Lie Algebras}

\author[Seok-Jin Kang, Jeong-Ah Kim and Dong-Uy Shin]
{Seok-Jin Kang$^{\star}$$^{*}$, Jeong-Ah Kim$^\dag$$^\diamond$ and
Dong-Uy Shin$^\dag$$^{*}$}

\address{$^*$School of Mathematics\\
         Korea Institute for Advanced Study\\
         Seoul 130-012, Korea}
\email{sjkang@kias.re.kr\\ shindong@kias.re.kr}
\address{$^\diamond$Department of Mathematics\\
         Seoul National University\\
         Seoul 151-747, Korea}
\email{jakim@math.snu.ac.kr}
\thanks{$^{\star}$This research was supported by KOSEF Grant
\# 98-0701-01-5-L and the Young Scientist Award, Korean Academy of
Science and Technology}
\thanks{$^{\dag}$This research was supported by KOSEF Grant
\# 98-0701-01-5-L and BK21 Mathematical Sciences Division, Seoul
National University}

\begin{abstract}
We give a new realization of crystal bases for finite dimensional
irreducible modules over special linear Lie algebras using the
monomials introduced by H. Nakajima. We also discuss the
connection between this monomial realization and the tableau
realization given by Kashiwara and Nakashima.
\end{abstract}
\maketitle
\vskip 1cm

\section*{Introduction}

The {\it quantum groups}, which are certain deformations of the
universal enveloping algebras of Kac-Moody algebras, were
introduced independently by V. G. Drinfeld and M. Jimbo
\cite{Drin,Jim}. In \cite{Kas90,Kas91}, M. Kashiwara developed the
{\it crystal basis theory} for integrable modules over quantum
groups. Crystal bases can be viewed as bases at $q=0$ and they are
given a structure of colored oriented graphs, called the {\it
crystal graphs}. Crystal graphs have many nice combinatorial
properties reflecting the internal structure of integrable
modules. Moreover, crystal bases have a remarkably nice behavior
with respect to taking the tensor product.

\vskip 2mm In \cite{Nakaj1}, while studying the structure of
quiver varieties, H. Nakajima discovered that one can define a
crystal structure on the set of irreducible components of a
lagrangian subvariety $\frak Z$ of the quiver variety $\frak M$.
These irreducible components are identified with certain
monomials, and the action of Kashiwara operators can be
interpreted as multiplication by monomials. Moreover, in
\cite{Nakaj1} and \cite{Kas02}, M. Kashiwara and H. Nakajima gave
a crystal structure on the set ${\mathcal M}$ of monomials and
they showed that the connected component ${\mathcal M}(\lambda)$
of ${\mathcal M}$ containing a highest weight vector $M$ with a
dominant integral weight $\lambda$ is isomorphic to the
irreducible highest weight crystal $B(\la)$. Therefore, a natural
question arises: for each dominant integral weight $\lambda$, can
we give an explicit characterization of the monomials in
${\mathcal M}(\lambda)$?

\vskip 2mm In this paper, for any dominant integral weight
$\lambda$, we give an explicit description of the crystal
${\mathcal M}(\lambda)$ for special linear Lie algebras. In
addition, we discuss the connection between the monomial
realization and tableau realization of crystal bases given by
Kashiwara and Nakashima. More precisely, let $T(\lambda)$ denote
the crystal consisting of semistandard tableaux of shape
$\lambda$. Then we show that there exists a canonical crystal
isomorphism between ${\mathcal M}(\lambda)$ and $T(\lambda)$,
which has a very natural interpretation in the language of
insertion scheme.

\vskip 2mm This work was initiated when the authors visited RIMS,
Kyoto University, in the spring of 2002. We would like to express
our sincere gratitude to Professor M. Kashiwara for his kindness
and stimulating discussions during our visit.


\vskip 1cm
\section{Crystal bases}

Let $I$ be a finite set and set $A=(a_{ij})_{i,j \in I}$ be a
generalized Cartan matrix. Consider the {\it  Cartan datum} $(A,
\Pi, \Pi^{\vee}, P, P^{\vee})$, where
$$\aligned
&P^{\vee}=(\bigoplus_{i\in I}\Z
h_i)\bigoplus(\bigoplus_{j=1}^{\text{corank} A}\Z d_j): \text{the
{\it dual weight lattice}},\\
&P=\{\lambda\in \mathfrak h^* | \lambda(P^{\vee})\subset \Z \}:
\text{the {\it weight lattice}},\\
&\Pi^{\vee} = \{ h_i | \, i\in I\}: \text{the set of {\it simple
coroots}},\\
&\Pi=\{\alpha_i | \, i\in I\}\subset \mathfrak h^*:  \text{the
{\it simple roots}}. \endaligned
$$
Let ${\mathfrak h}= \Q \otimes_{\Z} \pv$ be the {\it Cartan
subalgebra} and fix a nondegenerate symmetric bilinear form $( \ |
\ )$ on ${\mathfrak h}^*$ satisfying:
$$\frac{(\alpha_i | \alpha_i)}{2} \in \Z_{>0} \ \
\text{and} \ \ \lambda(h_i) = \frac{2
(\lambda|\alpha_i)}{(\alpha_i|\alpha_i)} \ \ \text{for all} \ i\in
I, \ \lambda \in {\mathfrak h}^*.$$

\vskip 2mm

The {\it quantum group} $U_q({\mathfrak g})$ associated with the
cartan datum  $(A, \Pi, \Pi^{\vee}, P, P^{\vee})$ is the
associative algebra over $\Q(q)$ with $1$ generated by the
elements
 $e_i$, $f_i$ $ (i \in I)$ and $q^h$ $(h \in P^{\vee})$ with the
following defining relations:
\begin{equation}
\begin{aligned}
\ & q^0 =1, \ \ q^h q^{h'} = q^{h+ h'} \quad \text{for}\,\, h, h'\in P^{\vee},\\
\ & q^h e_i q^{-h} = q^{\ali(h)} e_i, \quad
q^h f_i q^{-h} = q^{-\ali(h)} f_i \quad \text{for}\,\,h\in P^{\vee}, i\in I, \\
\ & e_i f_j - f_j e_i = \delta_{ij} \frac{K_i - K_i^{-1}}{q_i -
q_i^{-1}}
\quad \text{for}\,\,i,j \in I, \\
\ & \sum_{k=0}^{1-a_{ij}} (-1)^k {\begin{bmatrix} 1-a_{ij} \\ k
\end{bmatrix}}_{i}
e_i^{1-a_{ij}-k} e_j e_i^{k} = 0 \quad \text{for}\,\,i \neq j, \\
\ & \sum_{k=0}^{1-a_{ij}} (-1)^k {\begin{bmatrix} 1-a_{ij} \\ k
\end{bmatrix}}_{i} f_i^{1-a_{ij}-k} f_j f_i^{k} = 0 \quad \text{for}\,\,i
\neq j.
\end{aligned}
\end{equation}
Here, we use the notations:
\begin{equation*}
\begin{aligned}
& q_i = q^{\frac{(\alpha_i|\alpha_i)}{2}}, \quad K_i=
q^{\frac{(\alpha_i|\alpha_i)}{2}h_i}, \\
& [k]_i = \frac{q_i^k - q_i^{-k}}{q_i - q_i^{-1}}, \quad [n]_i !
=\prod_{k=1}^{n} [k]_i, \quad
{\begin{bmatrix} m \\ n \end{bmatrix}}_{i}
=\frac{[m]_{i}!}{[n]_{i}! \, [m-n]_{i}!}.
\end{aligned}
\end{equation*}
We also define $Q=\bigoplus_{i\in I}\Z\alpha_i$, $Q_+=\sum_{i\in
I}\Z_{\geq 0}\alpha_i$ and
\begin{center}
$P^{+} = \{ \lambda \in P \mid \lambda(h_i) \ge 0 \ \ \text{for
all} \ i\in I \}.$
\end{center}
In particular, the linear
functional $\Lambda_i\in P^{+}$ ($i\in I$) defined by
\begin{equation}
\Lambda_i(h_j)=\delta_{ij},\,\,\,\Lambda_i(d_s)=0\,\,\,\text{for}\,\,\,j\in
I, s=1,\cdots,\text{corank}A
\end{equation}
is called the {\it fundamental weight}.

\vskip 2mm The {\it category ${\mathcal O}_{int}$} consists of
$U_q(\g)$-modules $M$ satisfying the properties:
\begin{itemize}
\item [(i)] $M = \bigoplus_{\lambda \in P} M_{\lambda}$ with dim$M_{\lambda}<\infty$, where
$$
M_{\la} = \{ v\in M \mid q^h v = q^{\la(h)} v \text{ for all }
h\in \pv \},$$

\item [(ii)] there exist finitely many elements
$\lambda_1,\cdots, \lambda_s\in P$ such that
\begin{center}
$\text{wt}(M)\subset \bigcup_{j=1}^{s}(\lambda_j-Q_+),$
\end{center}
where $\text{wt}(M)=\{\la\in P| M_{\la}\neq 0\}$,
\item [(iii)]$e_i$ and $f_i$ ($i\in I$) are locally nilpotent
on $M$.
\end{itemize}
\noindent For each $i\in I$, it is well-known that every
$U_q(\g)$-module in the category ${\mathcal O}_{int}$ is a direct
sum of finite dimensional irreducible $U_{(i)}$-submodules, where
$U_{(i)}=\langle e_i,f_i,K_i^{\pm1}\rangle \cong
U_q(\mathfrak{sl}_2)$.

\vskip 3mm Fix an index $i\in I$ and set $e_i^{(n)} =
e_i^n/[n]_i!$, $f_i^{(n)} = f_i^n/[n]_i!$. Then every weight
vector $v\in M_{\lambda}$ can be written uniquely as
\begin{equation*}
v = \sum_{k\geq0} f_i^{(k)} v_k,
\end{equation*}
with $v_k \in \ker e_i \cap M_{\la+k\ali}$. We define the {\it
Kashiwara operators} $\eit$ and $\fit$ on $M$ by
\begin{equation}
\eit v = \sum_{k\geq1} f_i^{(k-1)} v_k, \qquad \fit v =
\sum_{k\geq0} f_i^{(k+1)} v_k.
\end{equation}

Let $\A_0 = \{ f/g \in \Q(q)  \, | \, f, g \in \Q[q], \, g(0)\neq
0 \}$.

\begin{defi}
A {\it crystal basis} of $M$ is a pair $(L, B)$ satisfying the
following conditions:

(i) $L$ is a free $\A_0$-submodule of $M$ such that $M \cong \Q(q)
\otimes_{\A_0} L$,

(ii) $B$ is a $\Q$-basis of $L/qL \cong \Q \otimes_{\A_0} L$,

(iii) $L=\bigoplus_{\lambda \in P} L_{\la}$, where $L_{\la} = L
\cap M_{\la}$,

(iv) $B=\bigsqcup_{\la \in P} B_{\la}$, where $B_{\la}=B \cap
\left(L_{\la} / q L_{\la} \right)$,

(v) $\eit L \subset L$, $\fit L \subset L$ for all $i\in I$,

(vi) $\eit B\subset B\cup \{0\}$, \ $\fit B\subset B\cup \{0\}$
for all $i\in I$,

(vii) for all $b, b'\in B$ and $i\in I$, $\fit b = b'$ if and only
if $b = \eit b'$.

\end{defi}

\noindent The set $B$ is given a colored oriented graph structure
with the arrows defined by
\begin{center}
$b \stackrel{i} \longrightarrow b' \quad \text{if and only if}
\quad \fit b = b'.$
\end{center}
The graph $B$ is called the {\it crystal graph} of $M$ and it
reflects the combinatorial structure of $M$. For instance, we have
$$ \dim_{\Q(q)} M_{\lambda} =\# B_{\lambda} \quad \text{for all
$\lambda \in P$}.$$
Moreover, the crystal basis have a very nice
behavior with respect to the tensor product. For each $b\in B$ and
$i\in I$, we define
\begin{equation}
\varepsilon_i(b)=\max \{k\ge 0 | \, \eit^k b \in B \}, \qquad
\varphi_i(b) = \max \{ k\ge 0 | \, \fit^k b \in B \}.
\end{equation}
Then we have:

\begin{prop} {\rm \cite{Kas90, Kas91}} \,\,
Let $M_j$ $(j=1, 2)$ be a  $U_q(\g)$-module in the category
${\mathcal O}_{int}$ and $(L_j, B_j)$ be its crystal basis. Set
\begin{equation*}
L = L_1 \otimes_{\A_0} L_2, \quad B= B_1 \times B_2.
\end{equation*}
Then $(L,B)$ is a crystal basis of $M_1 \otimes_{\Q(q)} M_2$,
where the Kashiwara operators on $B$ are given by
\begin{equation*}
\begin{aligned}
\eit(b_1\ot b_2)
  &=  \begin{cases}
  \eit b_1 \ot b_2 & \text{if $\vphi_i(b_1) \geq \vep_i(b_2)$,}\\
  b_1\ot\eit b_2 & \text{if $\vphi_i(b_1) < \vep_i(b_2)$,}
  \end{cases}\\
\fit(b_1\ot b_2)
  &=  \begin{cases}
  \fit b_1 \ot b_2 & \text{if $\vphi_i(b_1) > \vep_i(b_2)$,}\\
  b_1\ot\fit b_2 & \text{if $\vphi_i(b_1) \leq \vep_i(b_2)$.}
  \end{cases}
\end{aligned}
\end{equation*}
\end{prop}

\vskip 3mm We close this section with the existence and uniqueness
theorem for crystal bases.

\begin{prop} {\rm \cite{Kas91}} \label{prop:existence}\,\,
Let $V(\la)$ be the irreducible highest weight $U_q(\g)$-module
with highest weight $\la \in P^{+}$ and highest weight vector
$v_{\la}$. Let $L(\la)$ be the free $\A_0$-submodule of $V(\la)$
spanned by the vectors of the form $\tilde f_{i_1} \cdots \tilde
f_{i_r} v_{\la}$ $(i_k \in I, r\in \Z_{\ge 0})$ and set
\begin{equation*}
B(\la) = \{ \tilde f_{i_1} \cdots \tilde f_{i_r} u_{\la}+q L(\la)
\in L(\la) / q L(\la) \} \setminus \{0\}.
\end{equation*}
Then $(L(\la), B(\la))$ is a crystal basis of $V(\la)$ and every
crystal basis of $V(\la)$ is isomorphic to $(L(\la), B(\la))$.
\end{prop}

\vskip 1cm
\section{Nakajama's monomials}

\vskip 3mm  In this section, we briefly recall the crystal
structure on the set of monomials discovered by H. Nakajima
\cite{Nakaj1}. Our expression follows that of M. Kashiwara
\cite{Kas02}.

\vskip 2mm Let $\mathcal M$ be the set of monomials in the
variables $Y_i(n)$ for $i\in I$ and $n\in {\bf Z}$. Here, a
typical elements $M$ of ${\mathcal M}$ has the form
\begin{equation}
M=Y_{i_1}(n_1)^{a_1}\cdots Y_{i_r}(n_r)^{a_r},
\end{equation}
where $i_k\in I, n_k, a_k \in {\bf Z}$ for $k=1,\cdots, r.$ Since
$Y_i(n)$'s are commuting variables, we may assume that $n_1\leq
n_2\leq\cdots \leq n_r.$

For a monomial $M=Y_{i_1}(n_1)^{a_1}\cdots Y_{i_r}(n_r)^{a_r}$, we
define
\begin{equation} \aligned
\text{wt}(M)&=\sum_{k=1}^{r}a_k\La_{i_k}=a_1\La_{i_1}+\cdots a_r\La_{i_r},\\
\varphi_i(M)&=\text{max}\,\big(\big\{\sum_{k=1 \atop i_k=i}^s a_k \,\,|\,\,1\leq s\leq r \big\}\cup \{0\}\big),\\
\varepsilon_i(M)&=\text{max}\,\big(\big\{-\sum_{k=s+1 \atop
i_k=i}^r a_k \,\,|\,\,1\leq s\leq r-1 \big\}\cup \{0\}\big).
\endaligned
\end{equation}
It is easy to verify that $\varphi_i(M)\geq 0,\varepsilon_i(M)\geq
0,$ and $\langle h_i,
\text{wt}M\rangle=\varphi_i(M)-\varepsilon_i(M)$.

\vskip 2mm First, we define
\begin{equation}
\aligned
n_f&=\text{smallest}\,\,n_s\,\, \text{such that}\,\, \varphi_i(M)=\sum_{k=1 \atop i_k=i}^s a_k\\
&=\text{smallest}\,\,n_s\,\, \text{such that}\,\,
\varepsilon_i(M)=-\sum_{k=s+1 \atop i_k=i}^r a_k,\\
n_e&=\text{largest}\,\,n_s\,\, \text{such that}\,\, \varphi_i(M)=\sum_{k=1 \atop i_k=i}^s a_k\\
&=\text{largest}\,\,n_s\,\, \text{such that}\,\,
\varepsilon_i(M)=-\sum_{k=s+1 \atop i_k=i}^r a_k.
\endaligned
\end{equation}
In addition, choose a set $C=(c_{ij})_{i\neq j}$ of integers such
that $c_{ij}+c_{ji}=1,$ and define
$$A_i(n)=Y_i(n)Y_i(n+1)\prod_{j\neq
i}Y_j(n+c_{ji})^{\alpha_i(h_j)}.$$

Now, the {\it Kashiwara operators} $\eit$, $\fit$ ($i\in I$) on
${\mathcal M}$ are defined as follows:
\begin{equation}
\aligned
&\fit(M)=
  \begin{cases}
    0 & \text{if $\varphi_i(M)=0$}, \\
    A_i(n_f)^{-1}M & \text{if $\varphi_i(M)>0$},
  \end{cases}\\
&\eit(M)=
  \begin{cases}
    0 & \text{if $\varepsilon_i(M)=0$}, \\
    A_i(n_e)M & \text{if $\varepsilon_i(M)>0$}.
  \end{cases}
\endaligned
\end{equation}
Then the maps $\text{wt}: {\mathcal M}\rightarrow P$, $\varphi_i,
\varepsilon_i:{\mathcal M}\rightarrow \Z\cup\{-\infty\},$
$\eit,\fit:{\mathcal M}\rightarrow {\mathcal M}\cup\{0\}$ define a
$U_q(\mathfrak g)$-crystal structure on ${\mathcal M}$
\cite{Kas02,Nakaj1}.

\vskip 3mm Moreover, we have

\begin{prop}\, \cite{Kas02} \,\, \label{prop:mono}

{\rm (i)} For each $i\in I$, ${\mathcal M}$ is isomorphic to a
crystal graph of an integrable $U_{(i)}$-module.

{\rm (ii)} Let $M$ be a monomial with weight $\la$ such that $\eit
M=0$ for all $i\in I$, and let ${\mathcal M}(\la)$ be the
connected component of ${\mathcal M}$ containing $M$. Then there
exists a crystal isomorphism  $${\mathcal M}(\la) {\stackrel
{\sim}{\longrightarrow}} B(\la) \,\,\,\text{given by}\,\,\,
M\longmapsto v_{\la}.$$
\end{prop}

\vskip 2mm
\begin{example}
Let $\mathfrak g=A_2,$ and choose $c_{12}=1$ and $c_{21}=0.$ The
crystal ${\mathcal M}(\la)$ is given as follows.

(1)\,\,${\mathcal M}(\La_1):$
$$ Y_1(0) {\stackrel {1}{\longrightarrow}}Y_1(1)^{-1}Y_2(0){\stackrel {2}{\longrightarrow}}Y_2(1)^{-1}$$

\vskip 5mm (2)\,\, ${\mathcal M}(2\La_1):$
$$Y_1(0)^2 {\stackrel {1}{\longrightarrow}}Y_1(0)Y_1(1)^{-1}Y_2(0){\stackrel
{2}{\longrightarrow}}Y_1(0)Y_2(1)^{-1}$$ $$\hskip 1.5cm\downarrow
1 \hskip 3cm \downarrow 1 $$ $$\hskip 3.5cm Y_1(1)^{-2}Y_2(0)^2
\,\,{\stackrel
{2}{\longrightarrow}}\,\,Y_1(1)^{-1}Y_2(0)Y_2(1)^{-1}\,\,{\stackrel
{2}{\longrightarrow}}\,\,Y_2(1)^{-2}$$

\vskip 5mm (3) \,\, ${\mathcal M}(\La_1+\La_2):$

\begin{center}
\begin{texdraw}
\fontsize{10}{10}\selectfont \drawdim em \setunitscale 1

\htext(0 0){$Y_1(0)Y_2(0)$}

\move(-1 -1)\avec(-2.5 -2.5)\htext(-2.5 -1.5){$1$}\move(1
-1)\avec(2.5 -2.5)\htext(2.5 -1.5){$2$}

\htext(-4 -3.5){$Y_1(1)^{-1}Y_2(0)^2$}\htext(5
-3.5){$Y_1(0)Y_1(1)Y_2(1)^{-1}$}

\move(-3.5 -5)\avec(-2 -6.5)\htext(-2 -5.5){$2$}\move(4
-5)\avec(2.5 -6.5)\htext(2.5 -5.5){$1$}

\htext(-3 -8){$Y_1(0)Y_2(0)^{-1}$}\htext(4
-8){$Y_2(0)Y_2(1)^{-1}$}

\move(-0.5 -9)\avec(1.5 -11)\htext(1.8 -10.5){$1$}\move(1
-9)\avec(-1 -11)\htext(-1.3 -10.5){$2$}

\htext(-4 -12){$Y_1(1)^{-1}Y_1(2)^{-1}Y_2(0)$}\htext(5
-12){$Y_1(1)Y_2(1)^{-2}$}

\move(-3 -13)\avec(-1.5 -14.5)\htext(-1.5 -13.5){$2$}\move(4
-13)\avec(2.5 -14.5)\htext(2.5 -13.5){$1$}

\htext(1 -15.5){$Y_1(2)Y_2(1)^{-1}$}

\end{texdraw}
\end{center}

\vskip 5mm Note that ${\mathcal M}(\La_1)\cong B(\La_1),$
${\mathcal M}(2\La_1)\cong B(2\La_1),$ and ${\mathcal
M}(\La_1+\La_2)\cong B(\La_1+\La_2),$ respectively.

\end{example}
\vskip 1cm
\section{Characterization of ${\mathcal M}(\la)$}

In this section, we give an explicit characterization of the
crystal ${\mathcal M}(\la)$ for special linear Lie algebras. Let
$I=\{1,\cdots ,n \}$ and let $A=(a_{ij})_{i,j\in I}$ be the
generalized Cartan matrix of type $A_n.$ Here, the entries of $A$
are given by
\begin{equation}
a_{ij} = \begin{cases}
2 \quad & \text{if} \ \ i=j, \\
-1 \quad & \text{if} \ \ |i-j|=1, \\
0 \quad & \text{otherwise.}
\end{cases}
\end{equation}
We define by $U_q(\mathfrak g)=U_q(\mathfrak{sl}_{n+1})$ the
corresponding quantum group. For simplicity, we take the set
$C=(c_{ij})_{i\neq j}$ to be
\begin{equation}
c_{ij} = \begin{cases}
0 \quad & \text{if} \ \ i>j, \\
1 \quad & \text{if} \ \ i<j, \\
\end{cases}
\end{equation}
and set $Y_0(m)^{\pm 1}=Y_{n+1}(m)^{\pm 1}=1$ for all $m\in \Z.$
Then for $i\in I$ and $m\in \Z,$ we have
\begin{equation}
A_i(m)=Y_i(m)Y_i(m+1)Y_{i-1}(m+1)^{-1}Y_{i+1}(m)^{-1}.
\end{equation}
To characterize ${\mathcal M}(\la)$, we first focus on the case
when $\la=\La_k.$ Let $M_0=Y_k(m)$ for $m\in \Z.$ By $(2.2)$, we
see that
$$\text{wt}(M_0)=\La_k,\,\, \varphi_i(M_0)=\delta_{ik}\,\,\text{and}\,\,
\varepsilon_i(M_0)=0\,\,\,\text{for all}\,\,i\in I.$$ Hence $\eit
M_0=0$ for all $i\in I$ and the connected component containing
$M_0$ is isomorphic to $B(\La_k)$ over $U_q(\mathfrak g).$ For
simplicity, we will take $M_0=Y_k(0)$, even if that does not make
much difference.

\begin{prop}
For $k=1,\cdots ,n,$ let $M_0=Y_k(0)$ be a highest weight vector
of weight $\La_k.$ Then the connected component ${\mathcal
M}(\La_k)$ of ${\mathcal M}$ containing $M_0$ is characterized as
\begin{equation*}
{\mathcal M}(\La_k)=\Biggl\{\,\, \prod_{j=1}^r
Y_{a_j}(m_{j-1})^{-1}Y_{b_j}(m_j)\mid \begin{array}{l} {\rm (i)}\,
0\leq a_1<b_1<a_2<\cdots<a_r<b_r\leq n+1,\\{\rm
(ii)}\,k=m_0>m_1>\cdots>m_{r-1}>m_r=0,\\
{\rm (iii)}\, a_j+m_{j-1}=b_j+m_j\,\,\text{for all}
\,\,j=1,\cdots,r\le k.
\end{array}
 \,\,\Biggr\}.
\end{equation*}
\end{prop}

\begin{proof} By Proposition \ref{prop:mono}, it suffices to prove the following
statements:
\begin{itemize}
\item [(a)] For all $i\in I,$ we have  $\tilde e_i {\mathcal M(\la)} \subset {\mathcal M}(\la) \cup \{0\},
\quad \tilde f_i {\mathcal M}(\la) \subset {\mathcal M}(\la) \cup
\{0\}.$
\item [(b)] For all $M \in {\mathcal M}(\la)$, there exist a sequence of indices
$i_1,\cdots, i_t$ in $I$ such that
\begin{center}
$\tilde{e}_{i_1}\cdots\tilde{e}_{i_t}M=M_0.$
\end{center}
\end{itemize}

\vskip 2mm Let $i\in I$ and $M=Y_{a_1}(m_0)^{-1}Y_{b_1}(m_1)\cdots
Y_{a_r}(m_{r-1})^{-1}Y_{b_r}(m_r)\in {\mathcal M}(\La_k)$. If
$i\neq b_j$ for all $j$, then $\varphi_i(M)=0$, which implies
$\fit M=0$.

If $i=b_j$ for some $j$, then $\varphi_i(M)=\varphi_{b_j}(M)=1,$
$n_f=b_j,$ and
$$A_{b_j}(m_j)=Y_{b_j-1}(m_j+1)^{-1}Y_{b_j}(m_j)Y_{b_j}(m_j+1)Y_{b_j+1}(m_j)^{-1}.$$
Hence we obtain
$$
\begin{aligned}
\fit M &=A_{b_j}(m_j)^{-1}M \\
&=Y_{a_1}(m_0)^{-1}Y_{b_1}(m_1)\cdots
Y_{a_j}(m_{j-1})^{-1}Y_{b_j-1}(m_j+1)\\
&\hskip 5mm \times
Y_{b_j}(m_j+1)^{-1}Y_{b_j+1}(m_j)Y_{a_{j+1}}(m_j)^{-1}Y_{b_{j+1}}(m_{j+1})\cdots
Y_{a_r}(m_{r-1})^{-1}Y_{b_r}(m_r).
\end{aligned}
$$
If $a_j<b_j-1$ and $b_j+1<a_{j+1}$, then since
$a_j+m_{j-1}=(b_j-1)+(m_j+1)$, $b_j+(m_j+1)=(b_j+1)+m_j$, it is
easy to see that $\fit M\in {\mathcal M}(\La_k).$ If $a_j=b_j-1,$
then since $a_j+m_{j-1}=(b_j-1)+(m_{j}+1)$, we have
$m_{j-1}=m_j+1$, which implies
$Y_{a_j}(m_{j-1})^{-1}Y_{b_j-1}(m_j+1)=1.$ On the other hand, if
$b_j+1=a_{j+1},$ then $Y_{b_j+1}(m_j)Y_{a_{j+1}}(m_j)^{-1}=1$ and
$b_j+(m_j+1)=a_{j+1}+m_j=b_{j+1}+m_{j+1}.$ Hence $\fit M\in
{\mathcal M}(\La_k).$

Similarly, one can prove $\eit M\in {\mathcal M}(\La_k)\cup
\{0\}.$

\vskip 2mm To prove (b), we have only to show that if $M\in
{\mathcal M}(\La_k)$ and $\eit M=0$ for all $i\in I$, then
$M=M_0=Y_k(0).$ But this is obvious, for otherwise, we would have
$\varepsilon_{a_j}(M)=1\neq 0.$
\end{proof}

\begin{remark}
If we take $M_0=Y_k(N),$ then we have only to modify the condition
for $m_j$'s as follows:
$$ k+N=m_0>m_1>\cdots >m_{r-1}>m_r=N .$$
\end{remark}

\vskip 2mm For $i\in I$ and $m\in \Z$, we introduce new variables
\begin{equation}
X_i(m)=Y_{i-1}(m+1)^{-1}Y_i(m).
\end{equation}
Using this notation, every monomial
$M=\prod_{j=1}^rY_{a_j}(m_{j-1})^{-1}Y_{b_j}(m_j)\in {\mathcal
M}(\La_k)$ may be written as
$$ M=\prod_{j=1}^rX_{a_j+1}(m_{j-1}-1)X_{a_j+2}(m_{j-1}-2)\cdots
X_{b_j}(m_j).$$ For example, we have
$M_0=Y_k(0)=X_1(k-1)X_2(k-2)\cdots X_k(0).$

\vskip2mm Now, it is straightforward to verify that we have
another characterization of the crystal ${\mathcal M}(\La_k).$

\begin{cor}
For $k=1,\cdots, n,$ we have
\begin{center}
${\mathcal M}(\La_k)= \{X_{i_1}(k-1)X_{i_2}(k-2)\cdots
X_{i_k}(0)\,|\,1\leq i_1<i_2<\cdots <i_k\leq n+1 \}.$
\end{center}
\end{cor}

\begin{remark}
If we take $M_0=Y_k(N),$ then we need to replace $X_i(m)$ by
$X_i(m+N).$ That is,
$$ {\mathcal M}(\La_k)= \{X_{i_1}(N+k-1)X_{i_2}(N+k-2)\cdots
X_{i_k}(N)\,|\,1\leq i_1<i_2<\cdots <i_k\leq n+1 \}.$$
\end{remark}

\vskip 2mm We now consider the general case.
\begin{defi}
Set $M=\prod_t Y_{a_t}(m_t)^{-1}Y_{b_t}(n_t)$ with
$a_t+m_t=b_t+n_t$.

(a) For each $k=0,\cdots,n-1$, we define $M(k)^+$ to be the
product of  $Y_{a_t}(m_t)^{-1}Y_{b_t}(n_t)$'s in $M$ with $n_t=k$;
that is,
$$
\begin{aligned}
M(k)^+&=\prod_{t: n_t=k} Y_{a_t}(m_t)^{-1}Y_{b_t}(n_t)\\
&=\prod_{t} Y_{a_t}(m_t)^{-1}Y_{b_t}(k).
\end{aligned}
$$

(b) For each $k=1,\cdots,n$, we define $M(k)^-$ to be the product
of $Y_{a_t}(m_t)^{-1}Y_{b_t}(n_t)$'s in $M$ with $m_t=k$; that is,
$$
\begin{aligned}
M(k)^-&=\prod_{t:m_t=k} Y_{a_t}(m_t)^{-1}Y_{b_t}(n_t)\\
&=\prod_{t} Y_{a_t}(k)^{-1}Y_{b_t}(n_t).
\end{aligned}
$$
\end{defi}

Now, for $M(k)^+=\prod_{t} Y_{a_t}(m_t)^{-1}Y_{b_t}(k)$, we denote
by $\la^+(M(k))$ the sequence $(b_{i_1},b_{i_2},\cdots,b_{i_r})$
whose terms are arranged in such a way that $n+1\ge b_{i_1}\ge
b_{i_2}\ge \cdots\ge b_{i_r}$. Similarly, for $M(k)^-=\prod_{t}
Y_{a_t}(k)^{-1}Y_{b_t}(n_t)$, we denote by $\la^-(M(k))$ the
sequence $(a_{j_1},a_{j_{2}},\cdots,a_{j_s})$ whose terms are
arranged in such a way that $n+1> a_{j_1}\ge a_{j_2}\ge\cdots\ge
a_{j_s}$.

\begin{defi}
Let $(\la_1,\cdots,\la_r)$ and $(\mu_1,\cdots,\mu_s)$ be the
sequences such that
\begin{center}
$\la_i\geq \la_{i+1}$ ($1\leq i \leq r-1$), \quad $\mu_j\geq
\mu_{j+1}$ ($1\leq j\leq s-1$).
\end{center}
We define $(\la_1,\cdots,\la_r)\prec(\mu_1,\cdots,\mu_s)$ if
\begin{center}
$r\leq s$ and $\la_i<\mu_i$ for all $i=1,\cdots,r$.
\end{center}

\end{defi}

\begin{example}
Let $M$ be a monomial
$(Y_1(2)^{-1}Y_3(0))\cdot(Y_0(2)^{-1}Y_1(1))\cdot(Y_2(1)^{-1}Y_3(0))$.
Then $$\aligned
&M(0)^+=(Y_1(2)^{-1}Y_3(0))\cdot(Y_2(1)^{-1}Y_3(0)),\\
&M(1)^+=Y_0(2)^{-1}Y_1(1),\\
&M(1)^-=Y_2(1)^{-1}Y_3(0),\\
&M(2)^-=(Y_1(2)^{-1}Y_3(0))\cdot(Y_0(2)^{-1}Y_1(1)).
\endaligned
$$
Moreover, the sequences $\la^+(M(0))=(3,3)$, $\la^+(M(1))=(1)$,
$\la^-(M(1))=(2)$, and $\la^-(M(2))=(1,0)$. Therefore,
$\la^+(M(1))\prec\la^-(M(1))$.
\end{example}

\begin{thm} \label{chla}
Let $\la=a_1\La_1+\cdots+a_n\La_n$ be a dominant integral weight
and let $M_0=Y_1(0)^{a_1}\cdots Y_n(0)^{a_n}$ be a highest weight
vector of weight $\la$ in ${\mathcal M}.$ The connected component
${\mathcal M}(\la)$ in ${\mathcal M}$ containing $M_0$ is
characterized as the set of monomials of the form $$\prod_i
Y_{a_t}(m_t)^{-1}Y_{b_t}(n_t)$$ with $a_t+m_t=b_t+n_t$ satisfying
the following conditions:
\begin{itemize}
\item[(i)] $\la^+(M(k))\prec \la^-(M(k))$ for $k=1,\cdots,n-1$.
\item[(ii)] If
$\la^+(M(k))=(b_{i_1},b_{i_2},\cdots,b_{i_r})$ and
$\la^-(M(k))=(a_{j_1},a_{j_{2}},\cdots,a_{j_s})$, then $s-r=a_k.$
\end{itemize}
\end{thm}

\begin{proof} As in Proposition 3.1, it suffices to prove the
following statements:

(a) For all $i=1, \cdots, n$, we have
$$\tilde e_i {\mathcal M(\la)} \subset {\mathcal M}(\la) \cup \{0\}, \quad
\tilde f_i {\mathcal M}(\la) \subset {\mathcal M}(\la) \cup
\{0\}.$$

(b) If $M\in {\mathcal M(\la)}$ and $\eit M=0$ for all $i\in I$,
then $M=M_0.$

\vskip 2mm We first prove the statement $(a)$. Let $i\in I$ and
$M=\prod_t Y_{a_t}(m_t)^{-1}Y_{b_t}(n_t)$ be a monomial of
${\mathcal M}(\la)$. Assume that $\tilde f_iM\neq 0$. Then $i=b_t$
for some $t$ and $\tilde{f}_{b_t}M$ is obtained from $M$ by
multiplying
$$A_{b_t}(n_t)^{-1}=Y_{b_t-1}(n_t+1)Y_{b_t}(n_t)^{-1}Y_{b_t}(n_t+1)^{-1}Y_{b_t+1}(n_t).$$
If we express as $M=Y_{a_t}(m_t)^{-1}Y_{b_t}(n_t)M'$, then
$\tilde{f}_{b_t}M$ is expressed as
$$\tilde{f}_{b_t}M=Y_{a_t}(m_t)^{-1}Y_{b_t-1}(n_t+1)Y_{b_t}(n_t+1)^{-1}Y_{b_t+1}(n_t)M'.$$
Note that $(\tilde{f}_{b_t}M)(k)^{-}=M(k)^{-}$ and
$(\tilde{f}_{b_{t}}M)(k)^{+}=M(k)^{+}$ unless $k=n_{t}$ and
$n_{t}+1$. At first, consider the case when $k=n_t.$ Let
$\la^-(M(n_t))=(a_{j_1},a_{j_2},\cdots ,a_{j_s})$ and
$\la^+(M(n_t))=(b_{i_1},\cdots,b_{i_p}=b_t,\cdots,b_{i_r})$. Since
$\la^-(M(n_t))> \la^+(M(n_t))$, $a_{j_p}>b_{i_p}=b_t$. If
$a_{j_p}>b_t+1$, then we have
$$\aligned
\la^-(\tilde
f_{b_t}M(n_t))&=(a_{j_1},\cdots,a_{j_p},\cdots,a_{j_s}),\\
\la^+(\tilde
f_{b_t}M(n_t))&=(b_{i_1},\cdots,b_{i_p}=b_t+1,\cdots,b_{i_r})
\endaligned
$$
If $a_{j_p}=b_t+1$, then we have
$$\aligned
\la^-(\tilde f_{b_t}M(n_t))&=(a_{j_1},\cdots,a_{j_{p-1}},
a_{j_{p+1}},\cdots,
a_{j_s}),\\
\la^+(\tilde
f_{b_t}M(n_t))&=(b_{i_1},\cdots,b_{j_{p-1}},b_{j_{p+1}},\cdots,b_{i_r})
\endaligned
$$
It is clear that $\tilde f_{b_t}M$ satisfies the condition (i) and
(ii). Secondly, for the case $k=n_t+1$, we have
$$
\la^-((\tilde f_{b_t}M)(n_t+1))=\la^-(M(n_t+1))\cup \{b_t\},\quad
\la^+((\tilde f_{b_t}M)(n_t+1))=\la^+(M(n_t+1))\cup \{b_t-1\}.
$$
Hence $\tilde f_{b_t}M\in {\mathcal M}(\la)$.

\vskip 2mm Similarly, we can prove that $\tilde e_i {\mathcal
M(\la)} \subset {\mathcal M}(\la) \cup \{0\}.$

\vskip 2mm To prove (b), suppose $M\in {\mathcal M(\la)}$ and
$\tilde e_i M=0$ for all $i\in I.$ Then by the definition of the
$\varepsilon_i(M)$, $M=\prod_t Y_0(b_t)^{-1}Y_{b_t}(0)$. Since
$\text{wt}(M)=a_1\La_1+\cdots +a_n\La_n$ and
$\text{wt}(Y_k(0))=\La_k$, we have $M=Y_1(0)^{a_1}\cdots
Y_n(0)^{a_n}=\prod_{t=1}^{a_1+\cdots+a_n}
Y_0(b_t)^{-1}Y_{b_t}(0).$
\end{proof}

\begin{remark}
The crystal ${\mathcal M}(\la)$ is obtained by multiplying
$a_k$-many monomials in ${\mathcal M}(\La_{k})$ ($k=1,\cdots,n$).
That is,
$${\mathcal M}(\la)=\{M=M_{1,1}\cdots M_{1,a_1}M_{2,1}\cdots M_{n,a_n}\mid M_{k,l}\in {\mathcal M}(\La_{k})\,\,\,
\text{for $1\le k\le n$, $1\le l\le a_k$}\}.$$
\end{remark}

\begin{example} \label{exam:la}
Let $\la$ be a dominant integral weight $\La_1+2\La_2+\La_3$ of
$A_4$ and let $M=Y_1(0)Y_1(1)Y_1(2)^{-1}Y_2(1)^{-1}Y_3(0)^3$. Then
$M$ can be expressed as
\begin{center}
$M=Y_0(3)^{-1}Y_3(0)Y_1(2)^{-1}Y_3(0)Y_0(2)^{-1}Y_1(1)Y_2(1)^{-1}Y_3(0)Y_0(1)^{-1}Y_1(0)$.
\end{center}
Therefore, we have
$$
\aligned
&M(0)^+=Y_0(3)^{-1}Y_3(0)Y_1(2)^{-1}Y_3(0)Y_2(1)^{-1}Y_3(0)Y_0(1)^{-1}Y_1(0),\\
&M(1)^+=Y_0(2)^{-1}Y_1(1),\quad M(2)^+=M(3)^+=1,\\
\endaligned
$$
and
$$
\aligned
&M(1)^-=Y_2(1)^{-1}Y_3(0)Y_0(1)^{-1}Y_1(0),\\
&M(2)^-=Y_1(2)^{-1}Y_3(0)Y_0(2)^{-1}Y_1(1),\\
&M(3)^-=Y_0(3)^{-1}Y_3(0),\\
&M(4)^-=1.\\
\endaligned
$$
It is easy to see that $M$ satisfies the conditions of Theorem
\ref{chla}. Therefore, $M\in {\mathcal M}(\la)$.
\end{example}

\begin{defi}
Set $M=\prod_j X_{b_j}(n_j)$.

(i) For each $k=1,\cdots, n-1$, we define $M(k)$ by the monomial
obtained by multiplying all $X_{b_j}(n_j)$ with $n_j=k$ in $M$,
that is,
$$M(k)=\prod_{j:n_j=k} X_{b_j}(n_j)=\prod_{j} X_{b_j}(k).$$

(ii) For $M(k)=\prod_{j} X_{b_j}(k)$, we define by $\la(M(k))$ the
sequence $(b_{j_1}, b_{j_2},\cdots,b_{j_s})$ whose terms are
arranged in such a way that $n+1\ge b_{j_1}\ge b_{j_2}\ge \cdots
\ge b_{j_s}$.

\end{defi}

\begin{cor} \label{cor}
Let $\la=a_1\La_1+\cdots+a_n\La_n$. Then ${\mathcal M}(\la)$ is
expressed as the set of monomials
$$M=\prod_{1\le i\le n+1 \atop 0\le j\le n-1} X_{i}(j)^{m_{ij}}$$
such that
\begin{itemize}
\item [(i)] for each $j=0,1,\cdots,n-1$,
$$\sum_{i=1}^{n+1}m_{ij}=a_{j+1}+\cdots+a_n,$$

\item [(ii)] for each $j=1,\cdots,n-1$,
$\la(M(j))\prec \la(M(j-1))$.
\end{itemize}
\end{cor}

\begin{example}
Let $\la$ be a dominant integral weight $\La_1+2\La_2+\La_3$ of
$A_4$ and let $M$ be a monomial
$Y_1(0)Y_1(1)Y_1(2)^{-1}Y_2(1)^{-1}Y_3(0)^3$ given in Example
\ref{exam:la}. Then $M$ can be expressed as
\begin{center}
$ X_{1}(2) X_{2}(1)^2 X_{1}(1) X_{3}(0)^3 X_{1}(0)$
\end{center}
and so
\begin{center}
$\sum_{i=1}^5 m_{i0}=4$, $\sum_{i=1}^5 m_{i1}=3$, $\sum_{i=1}^5
m_{i2}=1$ and $\sum_{i=1}^5 m_{i3}=0$.
\end{center}
Moreover, since
\begin{center}
$M(2)=X_{1}(2)$, $M(1)= X_{1}(1) X_{2}(1)^2$ and $M(0)= X_{1}(0)
X_{3}(0)^3$,
\end{center}
we know that $\la(M(j))\prec\la(M(j-1))$
for all $j=1,2,3$. Therefore, $M\in {\mathcal M}(\la)$.
\end{example}

\vskip 3mm Consider the condition (ii) in Corollary \ref{cor}. For
$M=\prod_{1\le i\le n+1 \atop 0\le j\le n-1} X_{i}(j)^{m_{ij}}$,
there are $m_{i,j}$ -many $i$ entries in the sequence $\la(M(j))$.
Therefore, the condition $\la(M(j))\prec \la(M(j-1))$ implies that
\begin{equation}
\aligned &m_{1,n}=0,\,\,\,m_{ij}=0\quad \text{for $2\leq i\leq
n+1$, $n-i+2\leq
j\leq n$},\\
&\sum_{k=i}^{n+1}m_{k,j}\leq \sum_{k=i+1}^{n+1}m_{k,j-1}\quad
\text{for $i=1,\cdots, n+1$, $j=1,\cdots, n$}.
\endaligned
\end{equation}

Therefore, Corollary \ref{cor} is expressed as follows:
\begin{cor}\label{mono:matrix}
Let $\la=a_1\La_1+\cdots+a_n\La_n$. Then ${\mathcal M}(\la)$ is
expressed as the set of monomials
$$M=\prod_{1\le i\le n+1 \atop 0\le j\le n-1} X_{i}(j)^{m_{ij}}$$
such that
\begin{itemize}
\item [(i)] $m_{1,n}=0,\,\,\,m_{ij}=0$\quad for $2\leq i\leq
n+1$, $n-i+2\leq j\leq n$,\\
\item [(ii)] ${\displaystyle \sum_{i=1}^{n+1}}m_{ij}=a_{j+1}+\cdots+a_n$\quad for each
$j=0,1,\cdots,n-1$,\\
\item [(iii)] ${\displaystyle \sum_{k=i}^{n+1}}m_{k,j}\leq {\displaystyle \sum_{k=i+1}^{n+1}}m_{k,j-1}$\quad
for $i=1,\cdots, n+1$, $j=1,\cdots, n$.
\end{itemize}
\end{cor}


\vskip 1.5cm
\section{The connection with Young tableaux}
In this section, we give the correspondence between monomial
realization and tableau realization of crystal base for the
classical Lie algebra $\frak g=A_n$. To prove the results in this
section, we will adopt the expression of monomials given in
Corollary \ref{cor}.

Before we give the correspondence between monomial realization and
tableau realization, we introduce certain tableaux with given
shape which is different from Young diagram given by Kashiwara and
Nakashima.

\begin{defi}
(i) We define a {\it reverse Young diagram} to be a collection of
boxes in right-justified rows with a weakly decreasing number of
boxes in each row from bottom to top.

(ii) We define a (reverse) {\it tableau} by a reverse Young
diagram filled with positive integers.

(iii) A (reverse) tableau  $S$ is called a (reverse) {\it
semistandard tableau} if the entries in $S$ are weakly increasing
from left right in each row and strictly increasing from top to
bottom in each column.
\end{defi}

\begin{remark} A reverse Young diagram is just a diagram obtained by
reflecting Young diagram to the origin.
\end{remark}

Let $\la$ be a dominant integral weight. Let $S(\la)$ (resp.
$T(\la)$) be the set of all (reverse) semistandard tableaux (resp.
semistandard tableaux) of shape $\la$ with entries on
$\{1,2,\cdots,n\}$, which is realized as crystal basis of finite
dimensional irreducible modules \cite{KasNak,KS2}. For the
fundamental weight $\La_k$ $(k=1,\cdots,n)$, we have
$T(\La_k)=S(\La_k)$.

\begin{thm}\label{anyla}
Let $\la=a_1\La_1+\cdots+a_n\La_n$ be a dominant integral weight.
Then there is a crystal isomorphism $\psi:{\mathcal
M}(\la)\rightarrow S(\la)$.
\end{thm}

\begin{proof} Let $M$ be a monomial in ${\mathcal M}(\la)$. Then $M$ is
expressed as
$$M=\prod_{1\le i\le n+1 \atop 0\le j\le n-1} X_{i}(j)^{m_{ij}}.$$ We
define $\psi(M)$ to be the semistandard tableau with $m_{ij}$-many
$i$ entries in $(j+1)$-st row (from bottom to top) for
$i=1,\cdots, n+1$, $j=0,1,\cdots,n-1$. Indeed, by the condition
(ii) of Corollary \ref{mono:matrix}, the tableau $\psi(M)$ is of
shape $\la$. Moreover, the condition (i) and (iii) imply that
$\psi(M)$ is semistandard.

Conversely, let $S$ be a tableau of $S(\la)$ with $m_{i,j}$-many
$i$ entries in the $j$-th row (from bottom to top) for
$i=1,\cdots,n+1$ and $j=1,\cdots,n$. We define $\psi^{-1}(S)$ by
the monomial
$$\prod_{1\le i\le n+1 \atop 1\le j\le n} X_i(j-1)^{m_{i,j}}.$$
Then since $S$ is semistandard, it is easy to see that
$\psi^{-1}(S)$ satisfies the condition (i)-(iii) of Corollary
\ref{mono:matrix}. Moreover, it is clear that $\psi$ and
$\psi^{-1}$ are inverses of each other.

Now, it remains to show that $\psi$ is a crystal morphism. Let
$M=\prod X_i(j)^{m_{ij}}$ be a monomial in ${\mathcal M}(\la)$.
Let $X_a(k_1)$ and $X_a(k_2)$ be the monomials corresponding to
the entries $a$ in $S_{i,j}$ and $S_{i',j'}$, respectively. By the
definitions of $\psi(M)$ and $S_{i,j}$ ($1\le i\le n$, $1\le j\le
a_k$), we have the following fact:
\begin{center}
If $k_1>k_2$, then $i>i'$, or  $i=i'$ and $j>j'$.
\end{center}
Therefore, from the definition of Kashiwara operators on the set
${\mathcal M}$ of monomials and the tensor product rule of
Kashiwara operators which is applied to the set $S(\la)$, it is
easy to see that $\psi$ is a crystal morphism of
$U_q(A_n)$-modules.
\end{proof}

\savebox{\tmpfiga}{\begin{texdraw} \fontsize{7}{7}\selectfont
\drawdim em \setunitscale 0.85

\rlvec(2 0)\rlvec(0 6)\rlvec(-2 0)\rlvec(0 -6)\move(0 2)\rlvec(2
0)\move(0 4)\rlvec(2 0)

\htext(1 1){$4$}\htext(1 3){$3$}\htext(1 5){$2$}
\end{texdraw}}
\savebox{\tmpfigb}{\begin{texdraw} \fontsize{7}{7}\selectfont
\drawdim em \setunitscale 0.85

\rlvec(2 0)\rlvec(0 4)\rlvec(-2 0)\rlvec(0 -4)\move(0 2)\rlvec(2
0)

\htext(1 1){$4$}\htext(1 3){$1$}
\end{texdraw}}
\savebox{\tmpfigc}{\begin{texdraw} \fontsize{7}{7}\selectfont
\drawdim em \setunitscale 0.85

\rlvec(2 0)\rlvec(0 4)\rlvec(-2 0)\rlvec(0 -4)\move(0 2)\rlvec(2
0)

\htext(1 1){$2$}\htext(1 3){$1$}
\end{texdraw}}
\savebox{\tmpfigd}{\begin{texdraw} \fontsize{7}{7}\selectfont
\drawdim em \setunitscale 0.85

\rlvec(2 0)\rlvec(0 2)\rlvec(-2 0)\rlvec(0 -2)

\htext(1 1){$2$}
\end{texdraw}}
\savebox{\tmpfige}{\begin{texdraw} \fontsize{7}{7}\selectfont
\drawdim em \setunitscale 0.85

\rlvec(2 0)\rlvec(0 6)\move(0 0)\rlvec(0 6)\rlvec(8 0)\rlvec(0
-2)\rlvec(-8 0)\move(0 2)\rlvec(6 0)\rlvec(0 4)\move(4 2)\rlvec(0
4)

\htext(1 1){$4$}\htext(1 3){$2$}\htext(1 5){$1$}\htext(3
3){$2$}\htext(3 5){$1$}\htext(5 3){$3$}\htext(5 5){$2$}\htext(7
5){$4$}
\end{texdraw}}
\savebox{\tmpfigf}{\begin{texdraw} \fontsize{7}{7}\selectfont
\drawdim em \setunitscale 0.85

\rlvec(8 0)\rlvec(0 6)\rlvec(-2 0)\rlvec(0 -6)\move(8 4)\rlvec(-6
0)\rlvec(0 -4)\move(8 2) \rlvec(-8 0)\rlvec(0 -2)\move(4
0)\rlvec(0 4)

\htext(1 1){$2$}\htext(3 1){$2$}\htext(5 1){$4$}\htext(7
1){$4$}\htext(3 3){$1$}\htext(5 3){$1$}\htext(7 3){$3$}\htext(7 5
){$2$}
\end{texdraw}}

\begin{example}\label{map:mono-matrix}
Let $\la$ be a dominant integral weight $\La_1+2\La_2+\La_3$ of
$A_3$ and let $M$ be a monomial $Y_1(3)^{-1}Y_2(0)^2Y_3(1)^{-1}$,
then it can be expressed as
\begin{center}
$M=Y_1(3)^{-1}Y_4(0)(Y_0(2)^{-1}Y_2(0))^2 Y_3(1)^{-1}Y_4(0)$
\end{center} and so it is a monomial of ${\mathcal
M}(\La_1+2\La_2+\La_3)$. Moreover, $M$ is also expressed as
$$\aligned
M&=X_2(2)X_3(1)X_1(1)^2X_4(0)^2X_2(0)^2\\
&=(X_2(2)X_3(1)X_4(0))\times(X_1(1)X_4(0))\times(X_1(1)X_2(0))\times
X_2(0).
\endaligned$$
Then we have the matrix $(m_{ij})$ ($1\leq i\leq 4, \,\, 0\leq
j\leq 3$) and  the semistandard tableau $S\in
S(\La_1+2\La_2+\La_3)$ as follows:

$$
(m_{ij})=\left(
\begin{array}{cccc}
0 & 2 & 0 & 0 \\
2 & 0 & 1 & 0 \\
0 & 1 & 0 & 0 \\
2 & 0 & 0 & 0
\end{array} \right)
\quad \text{and} \quad S=\raisebox{-0.4\height}{\usebox{\tmpfigf}}
$$

\end{example}

\vskip 2mm
We have the following proposition between $S(\la)$ and
$T(\la).$

\begin{prop}\label{ks}\, \cite{KS,KS2} \,\,
For a dominant integral weight $\la=a_1\La_1+\cdots+a_n\La_n$,
there is a crystal isomorphism $\varphi: S(\la)\rightarrow T(\la)$
for $U_q(A_n)$-module given by $$\varphi(S)=S_{n,1}\leftarrow
S_{n,2}\leftarrow\cdots\leftarrow S_{n,a_n}\leftarrow
S_{n-1,1}\leftarrow\cdots\leftarrow S_{1,a_1},$$ where $S_{i,j}\in
S(\La_i)$ is the column of $S$ of length $i$ $(1\le i\le n$, $1\le
j\le a_i)$ from right to left.
\end{prop}

\begin{cor} Let $\la=a_1\La_1+\cdots+a_n\La_n$ be a dominant integral weight.
There is a crystal isomorphism $\phi:{\mathcal M}(\la)\rightarrow
T(\la).$
\end{cor}
\begin{proof}
By Theorem \ref{anyla} and Proposition \ref{ks},
$\phi=\varphi\circ \psi$ is a crystal isomorphism.
\end{proof}

\savebox{\tmpfiga}{\begin{texdraw} \fontsize{7}{7}\selectfont
\drawdim em \setunitscale 0.85

\rlvec(2 0)\rlvec(0 6)\rlvec(-2 0)\rlvec(0 -6)\move(0 2)\rlvec(2
0)\move(0 4)\rlvec(2 0)

\htext(1 1){$4$}\htext(1 3){$3$}\htext(1 5){$2$}
\end{texdraw}}
\savebox{\tmpfigb}{\begin{texdraw} \fontsize{7}{7}\selectfont
\drawdim em \setunitscale 0.85

\rlvec(2 0)\rlvec(0 4)\rlvec(-2 0)\rlvec(0 -4)\move(0 2)\rlvec(2
0)

\htext(1 1){$4$}\htext(1 3){$1$}
\end{texdraw}}
\savebox{\tmpfigc}{\begin{texdraw} \fontsize{7}{7}\selectfont
\drawdim em \setunitscale 0.85

\rlvec(2 0)\rlvec(0 4)\rlvec(-2 0)\rlvec(0 -4)\move(0 2)\rlvec(2
0)

\htext(1 1){$2$}\htext(1 3){$1$}
\end{texdraw}}
\savebox{\tmpfigd}{\begin{texdraw} \fontsize{7}{7}\selectfont
\drawdim em \setunitscale 0.85

\rlvec(2 0)\rlvec(0 2)\rlvec(-2 0)\rlvec(0 -2)

\htext(1 1){$2$}
\end{texdraw}}
\savebox{\tmpfige}{\begin{texdraw} \fontsize{7}{7}\selectfont
\drawdim em \setunitscale 0.85

\rlvec(2 0)\rlvec(0 6)\move(0 0)\rlvec(0 6)\rlvec(8 0)\rlvec(0
-2)\rlvec(-8 0)\move(0 2)\rlvec(6 0)\rlvec(0 4)\move(4 2)\rlvec(0
4)

\htext(1 1){$4$}\htext(1 3){$2$}\htext(1 5){$1$}\htext(3
3){$2$}\htext(3 5){$1$}\htext(5 3){$3$}\htext(5 5){$2$}\htext(7
5){$4$}
\end{texdraw}}
\savebox{\tmpfigf}{\begin{texdraw} \fontsize{7}{7}\selectfont
\drawdim em \setunitscale 0.85

\rlvec(8 0)\rlvec(0 6)\rlvec(-2 0)\rlvec(0 -6)\move(8 4)\rlvec(-6
0)\rlvec(0 -4)\move(8 2) \rlvec(-8 0)\rlvec(0 -2)\move(4
0)\rlvec(0 4)

\htext(1 1){$2$}\htext(3 1){$2$}\htext(5 1){$4$}\htext(7
1){$4$}\htext(3 3){$1$}\htext(5 3){$1$}\htext(7 3){$3$}\htext(7 5
){$2$}
\end{texdraw}}

\begin{example}
Let $M$ be a monomial $Y_1(3)^{-1}Y_2(0)^2Y_3(1)^{-1}$ of $A_3$
given in Example 4.4. Then we have
$$
\aligned
\phi(M)&=S_{3,1}\leftarrow S_{2,2}\leftarrow S_{2,1}\leftarrow S_{1,1}\\
&=\raisebox{-0.4\height}{\usebox{\tmpfiga}}\leftarrow
\raisebox{-0.4\height}{\usebox{\tmpfigb}}\leftarrow
\raisebox{-0.4\height}{\usebox{\tmpfigc}}\leftarrow
\raisebox{-0.3\height}{\usebox{\tmpfigd}}\\
&=\raisebox{-0.4\height}{\usebox{\tmpfige}}\,.
\endaligned
$$

\vskip 5mm Conversely, let $T$ be a tableau of
$T(\La_1+2\La_2+\La_3)$
$$T=\raisebox{-0.4\height}{\usebox{\tmpfige}}\,.$$
By applying the reverse bumping rule to the entries from bottom to
top and from right to left, i.e., from the entry $4$ of the
rightmost column to the entry $1$ on top of the leftmost column,
we have the following sequence
\begin{center}
$(2,3,4,1,4,1,2,2)$.
\end{center}
Therefore, we have
\begin{center}
$S_{3,1}=\raisebox{-0.4\height}{\usebox{\tmpfiga}}$\,,\,\,
$S_{2,2}=\raisebox{-0.4\height}{\usebox{\tmpfigb}}$\,,,\,\,
$S_{2,1}=\raisebox{-0.4\height}{\usebox{\tmpfigc}}$\, and
$S_{1,1}=\raisebox{-0.3\height}{\usebox{\tmpfigd}}$\,,
\end{center}
and since
\begin{center}
$\psi^{-1}(S_{3,1})=X_2(2)X_3(1)X_4(0)$,
$\psi^{-1}(S_{2,2})=X_1(1)X_4(0)$,\\
$\psi^{-1}(S_{2,1})=X_1(1)X_2(0)$, $\psi^{-1}(S_{1,1})=X_2(0)$,
\end{center}
we have
$$\aligned
\varphi^{-1}(T)&=\psi^{-1}(S_{3,1})\psi^{-1}(S_{2,2})\psi^{-1}(S_{2,1})\psi^{-1}(S_{1,1})\\
&=Y_1(3)^{-1}Y_4(0)(Y_0(2)^{-1}Y_2(0))^2 Y_3(1)^{-1}Y_4(0)\\
&=Y_1(3)^{-1}Y_2(0)^2Y_3(1)^{-1}.
\endaligned$$

\end{example}

\vskip 15mm


\begin{thebibliography}{1}


\bibitem{Drin}
V. G. Drinfeld, \emph{Hopf algebras and the quantum Yang-Baxter
equation}, Soviet Math. Dokl. \textbf{32} (1985), 254--258.

\bibitem{Fulton}
W.~Fulton, \emph{Young Tableaux {\rm :} with applications to
representation theory and geometry}, Cambridge University Press,
1997.

\bibitem{HK}
J.~Hong, S.-J. Kang, \emph{Introduction to Quantum Groups and
Crystal Bases},  Graduate Studies in Mathematics \textbf{42},
Amer. Math. Soc., 2002.

\bibitem{Jim}
M.~Jimbo, \emph{A $q$-difference analogue of $U(\frak g)$ and the
Yang-Baxter equation}, Lett. Math. Phys. \textbf{10} (1985),
63--69.

\bibitem{Kac90}
V.~Kac, \emph{Infinite Dimensional Lie Algebras}, Cambridge
University Press, 3rd ed., 1990.


\bibitem{Kas90}
M.~Kashiwara, \emph{Crystalizing the $q$-analogue of universal
enveloping
  algebras}, Comm. Math. Phys. \textbf{133} (1990), 249--260.

\bibitem{Kas91}
\bysame, \emph{On crystal bases of the $q$-analogue of universal
enveloping algebras}, Duke Math. J. \textbf{63} (1991), 465--516.

\bibitem{Kas02}
\bysame, \emph{Realizations of crystals}, to appear in Contemp.
Math.


\bibitem{KasNak}
M.~Kashiwara, T.~Nakashima, \emph{Crystal graphs for
representations of the $q$-analogue of classical Lie algebras}, J.
Algebra \textbf{165} (1994), 295-345.

\bibitem{KS}
J.-A. Kim, D.-U. Shin, \emph{Insertion scheme for the crystal of
the classical Lie algebras}, submitted

\bibitem{KS2}
\bysame, \emph{Correspondence between Young walls and Young
tableaux realizations of crystal bases for the classical Lie
algebras}, submitted

\bibitem{Mac}
I. G. Macdonald, \emph{Symmetric Functions and Hall Polynomials},
Oxford University Press, Oxford, 2nd ed., 1995.

\bibitem{Nakaj1}
H.~Nakajima, \emph{Quiver varieties and tensor products}, Invent.
Math. \textbf{146} (2001), 399--449.

\bibitem{Nakaj2}
H.~Nakajima, \emph{$t$-analogs of $q$-characters of quantum affine
algebras of type $A_n$, $D_n$}, to appear in Contemp. Math.

\bibitem{Nak}
T.~Nakashima, \emph{Crystal base and a generalization of the
Littlewood-Richardson rule for classical Lie algebras}, Comm.
Math. Phys. \textbf{154} (1993), 215--243.



\end{thebibliography}

\providecommand{\bysame}{\leavevmode\hbox
to5em{\hrulefill}\thinspace}

\end{document}